\documentclass{amsart}

\usepackage{amsmath}
\usepackage{amscd}
\usepackage{amssymb}
 \usepackage{color}

\input xy
\xyoption{all}

\newcommand{\pd}{\partial}

\newcommand{\bZ}{{\mathbb Z}}

\newcommand{\cP}{{\mathcal P}}
\newcommand{\half}{\frac{1}{2}}

\newcommand{\vac}{|0\rangle}

\newcommand{\cor}[1]{\langle {#1} \rangle}
 \newcommand{\bp}{{\mathbf p}}
\newcommand{\bx}{{\mathbf x}}\newcommand{\by}{{\mathbf y}}
\DeclareMathOperator{\Aut}{Aut} 
\DeclareMathOperator{\Id}{Id}

\newtheorem{theorem}{Theorem}[section]
\newtheorem{theorem/definition}{Theorem/Definition}[section]

\newtheorem{proposition}{Proposition}[section]
\newtheorem{lemma}{Lemma}[section]

\theoremstyle{remark}

\theoremstyle{definition}

\newcommand{\bml}{\begin{multline}}
\newcommand{\eml}{\end{multline}}

\newcommand{\bag}{\begin{align}}
\newcommand{\egn}{ \end{align}}

\newcommand{\be}{\begin{equation}}
\newcommand{\ee}{\end{equation}}
\newcommand{\bea}{\begin{eqnarray}}
\newcommand{\eea}{\end{eqnarray}}
\newcommand{\ben}{\begin{eqnarray*}}
\newcommand{\een}{\end{eqnarray*}}
\newcommand{\bet}{\begin{equation}
\begin{split}}
\newcommand{\eet}{\end{split}
\end{equation}}

\definecolor{yellow}{rgb}{1,1,0}
\definecolor{orange}{rgb}{1,.7,0}
\definecolor{red}{rgb}{1,0,0} \definecolor{blue}{rgb}{0,0,1}
\definecolor{white}{rgb}{1,1,1}

\definecolor{A}{rgb}{.75,1,.75}

\begin{document}

\title
{A Proof of the Full Mari\~no-Vafa Conjecture}
\author{Jian Zhou}
\address{Department of Mathematical Sciences\\Tsinghua University\\Beijing, 100084, China}
\email{jzhou@math.tsinghua.edu.cn}

\begin{abstract}
We present a proof of the full Mari\~no-Vafa Conjecture that identifies
certain open string invariants of the resolved conifold with the Chern-Simons knot invariant of the unknot,
i.e. the quantum dimensions.
\end{abstract}

% \date{\today}

\maketitle

\section{Introduction}

The duality between topological string theory and Chern-Simons theory
is one of the important examples of string duality.
Mathematically, this predicts an amazing connection between Gromov-Witten theory
and link invariants.
Such predictions have inspired many mathematical results that verify this connections.
In this short note we will present a proof of the full Mari\~no-Vafa Conjecture that identifies
certain open string invariants of the resolved conifold with the Chern-Simons knot invariant of the unknot,
i.e. the quantum dimensions.

Witten \cite{Wit} proposed that the large N expansion of Chern-Simons link invariants \cite{Wit1}
can be identified with the genus expansion of open string invariants of the deformed conifold $T^*S^3$,
which counts holomorphic maps from Riemann surfaces boundaries to $T^*S^3$,
with $h$ boundary components mapped to the Lagrangian submanifold $S^3$.
Gopakumar-Vafa \cite{Gop-Vaf} proposed that after summing over the number $h$ of boundary components,
one gets from the large $N$ expansion of the $3$-manifold invariant
of $S^3$ the genus expansion of closed string invariants of the resolved conifold.
Ooguri and Vafa \cite{Oog-Vaf} extended this and proposed that the large $N$ expansion of link invariants
corresponds to open string invariants of the resolved conifold,
associated with some Lagrangian submanifolds corresponding to the link.
Mari\~no-Vafa \cite{Mar-Vaf} extended this further by considering framed knots.
Later,
the duality has been extended to other local Calabi-Yau geometries \cite{AMV} and a formalism called
the topological vertex \cite{AKMV} has been developed to compute open and closed string invariants
of local toric Calabi-Yau geometries,
in terms of quantities related to the link invariants of the Hopf link.
There are many related works in the physics literature.
More recently,
duality between topological strings
and link homology \cite{Kho-Roz} has been proposed \cite{GSV, GIKV}.

Many predictions in the physics literature has been made mathematically rigorous.
By comparing with the formal localization calculations performed for the resolved conifold by Katz and Liu \cite{Kat-Liu},
Mari\~no-Vafa \cite{Mar-Vaf} conjectured a closed formula for some Hodge integrals.
This formula for Hodge integrals has been proved \cite{LLZ1, Oko-Pan}
and generalized \cite{Zho1, LLZ2, LLLZ}.
Such formulas turns out to be important for the mathematical
calculations of open and closed string invariants of local Calabi-Yau geometries \cite{Zho2, LLLZ}.
Because mathematically the open string invariants is very difficult to define in symplectic geometry,
relative moduli spaces in algebraic geometry have been used to bypass this difficulty  \cite{Li-Son, LLLZ}.

Our proof of the full Mari\~no-Vafa starts with the computations of the open string invariants
of the resolved conifold by the mathematical theory of the topological vertex \cite{LLLZ}.
Such computations involve complicated summations over two partitions,
they are difficult to carry out directly.
We will use a simplifying trick discovered in an earlier work \cite{Zho4}.

\vspace{.1in}
{\em Acknowledgements}.
This research is partially supported by two NSFC grants (10425101 and 10631050)
and a 973 project grant NKBRPC (2006cB805905).

\section{Preliminaries}

In this section we express the quantum dimensions as specializations
of Schur functions.

\subsection{Symmetric functions}

The Newton functions $\{p_\mu\}$ and the Schur functions $\{s_\nu\}$ form additive bases of $\Lambda$.
They are related by the character values:
\begin{align} \label{eqn:NewtonSchur}
p_{\mu} & = \sum_{\nu} \chi_{\nu}(\mu) s_{\nu}, & s_{\nu} & =  \sum_\mu \frac{\chi_\nu(\mu)}{z_\mu} p_\mu.
\end{align}
Here $\mu = (\mu_1, \dots, \mu_l)$ is a partition,
\be
z_\mu = \prod_{i=1}^{l(\mu)} \mu_i \cdot |\Aut(\mu)|,
\ee
and $\chi_\mu$ denotes the character of the irreducible representation $S_{|\mu|}$, $|\mu| = \sum_i \mu_i$,
indexed by $\mu$,
and $\chi_\mu(\nu)$ denotes the value of $\chi_\mu$ on the conjugacy class of $S_{|\mu|}$ indexed by $\nu$.
The character values satisfy the following orthogonality relations:
\begin{align}
\sum_\nu \frac{1}{z_\nu} \chi_\mu(\nu) \chi_{\eta}(\nu) & = \delta_{\mu, \eta},
& \sum_{\mu} \chi_\mu(\xi)\chi_\mu(\phi) = \delta_{\xi, \phi} z_\xi. \label{eqn:Orth2}
\end{align}
Let $\bx = (x_1, x_2, \dots)$ and $\by = (y_1, y_2, \dots)$.
Then by (\ref{eqn:NewtonSchur}) and (\ref{eqn:Orth2}) one has
\be \label{eqn:NewtonSchur2}
\sum_\mu s_\mu(\bx)s_\mu(\by) = \sum_\nu \frac{1}{z_\nu} p_\nu(\bx)p_\nu(\by).
\ee
%% The right-hand side can be rewritten as follows:
%% \be
%% \sum_\nu \frac{1}{z_\nu} p_\nu(\bx)p_\nu(\by)
%% = \exp \sum_{n=1}^\infty \frac{1}{n} p_n(\bx) p_n(\by).
%% \ee
%% So we get:
%% \be
%% \exp \sum_{n=1}^\infty \frac{1}{n} p_n(\bx) p_n(\by) = \sum_\mu s_\mu(\bx)s_\mu(\by).
%% \ee

\subsection{Scalar product and vacuum expectation values}

There is a natural scalar product on $\Lambda$ such that:
\begin{align}
\langle p_\mu, p_\nu\rangle & = \delta_{\mu, \nu} z_{\mu}, &
\langle s_\mu, s_\nu\rangle & = \delta_{\mu, \nu}.
\end{align}
In particular,
\be
\langle p_\mu, s_\nu \rangle = \chi_\nu(\mu).
\ee
Following physical notations,
we will write $|\mu\rangle := s_\mu$,
and write the inner product $\langle A s_\mu, s_\nu\rangle$ for some linear operator $A: \Lambda \to \Lambda$
as $\langle \nu | A| \mu\rangle$.
In particular,
$$\cor{A}:=\langle0|A|0\rangle$$
is called the vacuum expectation value,
where $\vac = 1 \in \Lambda$.

\subsection{Creators and annihilators}
On $\Lambda$ one can introduce the following operators for nonzero integers $n$:
\be
\beta_n: \Lambda \to \Lambda, \qquad \beta_n (f) = \begin{cases}
p_{-n} \cdot f, & n < 0, \\
n \frac{\pd}{\pd p_n} f, & n > 0.
\end{cases}
\ee
These operators satisfy the following Heisenberg commutation relations:
\be
[\beta_m, \beta_n ] = m \delta_{m, -n} \Id_{\Lambda}.
\ee
For a partition $\mu = (\mu_1, \dots, \mu_l)$,
we write
\begin{align}
\beta_\mu & = \prod_{i=1}^l \beta_{\mu_i}, & \beta_{-\mu} & = \prod_{i=1}^l \beta_{-\mu_i}.
\end{align}
One clearly has
\be
\beta_n \vac = 0, \quad n > 0;
\qquad\qquad
\beta_{-\mu} \vac = p_{\mu}.
\ee
Hence the operators $\{\beta_n\}_{n>0}$ are called the annihilators,
and the operators $\{\beta_{-n}\}_{n > 0}$ are called the creators.

It is easy to see that
\be
\cor{ \beta_\mu \beta_{-\nu}} = \delta_{\mu, \nu} z_{\mu}.
\ee
This is called the Wick formula.
It follows from this identity that
\be
\cor{\exp (\sum_{n=1} \frac{a_n}{n} \beta_n) \cdot \exp ( \sum_{n=1}^\infty
\frac{b_n}{n} \beta_{-n} ) }
= \exp \sum_{n=1}^\infty \frac{a_nb_n}{n}.
\ee
We will understand $\Lambda$ as the space of symmetric functions in $\bx$.
Then we have
\be
\exp (\sum_{n=1}^\infty \frac{1}{n} p_n(\by) \beta_{-n}) \vac
= \sum_\nu \frac{1}{z_\nu} p_\nu(\by) \beta_{-\nu} \vac.
\ee
Now we apply (\ref{eqn:NewtonSchur2}) to get:
\be \label{eqn:NewtonSchur3}
\exp (\sum_{n=1}^\infty \frac{1}{n} p_n(\by) \beta_{-n}) \vac
= \sum_\mu s_\mu(\by) |\mu\rangle.
\ee

\subsection{Involution}

One can define an involution $\omega: \Lambda \to \Lambda$ by
\be
\omega(p_n) = (-1)^{n-1} p_n.
\ee
One has
\be
\omega(s_\mu) = s_{\mu^t}.
\ee
Apply $\omega$ on both sides of (\ref{eqn:NewtonSchur3}):
\be \label{eqn:NewtonSchur3t}
\exp (\sum_{n=1}^\infty \frac{1}{n} (-1)^{n-1} p_n(\by) \beta_{-n}) \vac
= \sum_\mu s_{\mu^t}(\by) |\mu\rangle.
\ee

\subsection{Cut-and-join operator}
Another useful operator on $\Lambda$ is the cut-and-join operator
\begin{align}
K: & = \half \sum_{i,j=1}^\infty \big(p_{i+j} ij\frac{\pd^2}{\pd p_i\pd p_j}
+ p_ip_j (i+j) \frac{\pd}{\pd p_{i+j}} \big) \\
& = \half \sum_{i,j =1}^\infty (\beta_{-(i+j)} \beta_i\beta_j + \beta_{-i}\beta_{-j} \beta_{i+j}). \nonumber
\end{align}
The Schur functions are eigenvectors of this operator:
\be
K |\mu\rangle = \half \kappa_\mu |\mu\rangle.
\ee

\subsection{Specializations of symmetric functions}

When $\by = q^{\rho}:=(q^{-1/2}, q^{-3/2}, \dots)$,
it is easy to see that
\be
p_n(q^\rho) = \frac{1}{q^{n/2} - q^{-n/2}} = \frac{1}{[n]}.
\ee
It is a very interesting fact that
\be
s_\mu(q^\rho) = \frac{q^{\kappa_\mu/4}}{\prod_{x \in \mu} [h(x)]},
\ee
where
$\kappa_\mu = \sum_{i=1}^{l(\mu)} \mu_i(\mu_i - 2 i + 1)$.
In other words,
when $p_n(\by) = \frac{1}{[n]}$ for all $n \geq 1$,
one has
$$s_\mu(\by) = \frac{q^{\kappa_\mu/4}}{\prod_{x \in \mu} [h(x)]}.$$
By (\ref{eqn:NewtonSchur3}),
\be
\exp \sum_{n=1}^\infty \frac{1}{n[n]} \beta_{-n}
= \sum_\mu \frac{q^{\kappa_\mu/4}}{\prod_{x \in \mu} [h(x)]} |\mu\rangle.
\ee

We will use another useful specialization.
When
\be
p_n(\by) = \frac{a^n-b^n}{1-q^n}
\ee
for $n \geq 1$,
one has
\be
s_\mu(\by) = q^{n(\mu)} \prod_{x\in \mu}  \frac{a - b q^{c(x)}}{1- q^{h(x)}}.
\ee
Here we represent a partition $\mu$ by its Young diagram,
$\mu^t$ is the partition corresponding to the transposed Young diagram of $\mu$.
Define
\be
n(\mu) = \sum_{i=1}^{l(\mu)} (i-1) \mu_i
= \sum_{i=1}^{l(\mu^t)} \binom{\mu^t_i}{2}.
\ee
For a box $x\in \mu$ at the $i$-row and $j$-th column,
its content and hook length are defined by:
\begin{align}
c(x) & = j - i,
& h(x) & = \mu_i + \mu^t_j -i-j + 1,
\end{align}
respectively.
The following identities hold:
\bea
&& \sum_{x\in \mu} h(x) = n(\mu) + n(\mu^t) + |\mu|, \\
&& \sum_{x\in \mu} c(x) = n(\mu^t) - n(\mu) = \half \kappa_\mu.
\eea
With these facts,
one can easily get the following

\begin{lemma}
If for $n \geq 1$ one has
\be \label{eqn:Pn}
p_n(\by) = \frac{1 - e^{-nt}}{[n]},
\ee
then one has:
\be \label{eqn:Smu}
s_\mu(\by) = e^{-|\mu|t/2} \prod_{x \in \mu} \frac{[c(x)]_{e^t } }{[h(x)]}.
\ee
where
\be
[n]_{e^t} = e^{nt/2} q^{n/2} - e^{-nt/2}q^{-n/2}.
\ee
\end{lemma}

\subsection{Relationship with the large N quantum dimension}

The colored large N HOMFLY polynomials are given by the quantum dimension
given as follows \cite[(5.4)]{Mar-Vaf}:
\be
\dim_q R_\mu = \prod_{1 \leq i < j \leq l(\mu)} \frac{[\mu_i-\mu_j+j-i]}{[j-i]}
\cdot \prod_{i=1}^{l(\mu)}
\frac{\prod_{j=1}^{\mu_i} [j-i]_{e^t}}{\prod_{j=1}^{\mu_i} [j-i+l(\mu)]}.
\ee
By \cite[\S 1.1, Example 1]{Mac},
\ben
&& \sum_{x\in \mu} t^{h(x)} + \sum_{1 \leq i < j \leq l(\mu)} t^{\mu_i-\mu_j + j-i} \\
& = & \sum_{i=1}^{l(\mu)} \sum_{j=1}^{\mu_i - i + l(\mu)} t^j
= \sum_{j=1}^{\mu_i} t^{j-i +l(\mu)} + \sum_{1 \leq i < j \leq l(\mu)} t^{j-i},
\een
therefore,
\be \label{eqn:QDim}
\dim_q R_\mu = \prod_{x \in \mu} \frac{[c(x)]_{e^{t}}}{[h(x)]}.
\ee

By combining (\ref{eqn:NewtonSchur3}), (\ref{eqn:NewtonSchur3t}), (\ref{eqn:Pn}), (\ref{eqn:Smu}) with (\ref{eqn:QDim}),
we get:

\begin{proposition}
The following identities hold:
\bea
&& \exp (\sum_{n=1}^\infty \frac{1-e^{-nt} }{n [n]} \beta_{-n} ) \vac
= \sum_\mu e^{-|\mu|t/2}  \cdot \dim_q R_\mu \cdot |\mu\rangle,  \label{eqn:QDim2} \\
&& \exp (\sum_{n=1}^\infty (-1)^{n-1} \frac{1-e^{-nt} }{n [n]} \beta_{-n} ) \vac
= \sum_\mu e^{-|\mu|t/2} \cdot \dim_q R_{\mu^t} \cdot |\mu\rangle. \label{eqn:QDimT}
\eea
\end{proposition}

\section{Proof of the Full Mari\~no-Vafa Conjecture}

In \S 3.1 we first recall some results in \cite{Zho4}.
Then we prove the full Mari\~no-Vafa Conjecture in \S 3.2.

\subsection{Open string amplitudes with one outer brane}

We compute the open string amplitude for the resolved conifold with one outer brane and framing $a$
by the theory of the topological vertex \cite{AKMV, LLLZ}.
There are two possibilities,
corresponding to two different toric diagrams as follows:
$$
\xy
(0,10); (0,0), **@{-}; (-7, -7), **@{-}; (-10, -4), **@{.};
(0,0); (20,0), **@{-}; (27,7), **@{-}; (20,0); (20, -10), **@{-};
(-5, -2)*+{\mu}; (5,1.5)*+{\nu^t}; (15,1)*+{\nu}; (11,-14)*+{(a)};
\endxy
\qquad
\xy
(-4,5); (0,10), **@{.}; (0,0), **@{-}; (-7, -7), **@{-};
(0,0); (20,0), **@{-}; (27,7), **@{-}; (20,0); (20, -10), **@{-};
(-2, 3)*+{\mu}; (5,1.5)*+{\nu^t}; (15,1)*+{\nu}; (11,-14)*+{(b)};
\endxy
$$
\begin{center}
{\bf Figure 1.}
\end{center}
For Figure 1(a) we have
\be
Z^{(a)}(\lambda; t; \bp)
= \sum_{\mu, \nu, \eta} q^{a \kappa_\mu/2}  C_{(0), \mu, \nu^t}(q) \cdot Q^{|\nu|} \cdot
C_{\nu,(0), (0)}(q) \cdot
\frac{\chi_{\mu}(\eta)}{z_{\eta}\sqrt{-1}^{l(\eta)}} p_{\eta}(\bx),
\ee
where $q = e^{\sqrt{-1}\lambda}$ and $Q = - e^{-t}$.
For Figure 1(b) we have
\be
\tilde{Z}^{(a)}(\lambda; t; \bp)
= \sum_{\mu, \nu, \eta} q^{a \kappa_\mu/2}  C_{\mu, (0), \nu^t}(q) \cdot Q^{|\nu|} \cdot
C_{\nu,(0), (0)}(q) \cdot
\frac{\chi_{\mu}(\eta)}{z_{\eta}\sqrt{-1}^{l(\eta)}} p_{\eta}(\bx).
\ee
The topological vertices here can be rewritten as follows (cf \cite{Zho3}):
\be
C_{(0), \mu, \nu^t}(q) = q^{-\kappa_\nu/2} W_{\mu,\nu}(q),
\ee
so we have
\be \label{eqn:Za}
Z^{(a)}(\lambda;t;\bp) = \sum_{\mu, \nu, \eta} q^{a \kappa_\mu/2} W_{\mu\nu}(q) q^{-\kappa_\nu/2} Q^{|\nu|}  W_{\nu}(q)\cdot
\frac{\chi_{\mu}(\eta)}{z_{\eta}\sqrt{-1}^{l(\eta)}} p_{\eta}(\bx),
\ee
and
\be \label{eqn:Zb}
\tilde{Z}^{(a)}(\lambda;t;\bp) = \sum_{\mu, \nu, \eta} q^{(a+1) \kappa_\mu/2} W_{\mu^t\nu^t}(q) Q^{|\nu|}  W_{\nu}(q)\cdot
\frac{\chi_{\mu}(\eta)}{z_{\eta}\sqrt{-1}^{l(\eta)}} p_{\eta}(\bx).
\ee
The normalized open string amplitudes are defined by:
\begin{align}
\hat{Z}^{(a)}(\lambda;t;\bp) = \frac{Z^{(a)}(\lambda;t; \bp)}{Z^{(a)}(\lambda;t; \bp)|_{\bp = {\bf 0}}}, &&
\hat{\tilde{Z}}^{(a)}(\lambda;t;\bp) = \frac{\tilde{Z}^{(a)}(\lambda;t; \bp)}{\tilde{Z}^{(a)}(\lambda;t; \bp)|_{\bp = {\bf 0}}}.
\end{align}
In \cite{Zho4} we have proved the following:

\begin{proposition}
The normalized open string amplitude of the resolved conifold with one outer brane and framing $a \in \bZ$
is given by:
\be \label{eqn:ZaOperator1}
\hat{Z}^{(a)}(\lambda;t; \bp)
=  \langle \exp (\sum_{n=1}^\infty \frac{p_n(\bx)}{ni} \beta_n) q^{(a+1)K}
\exp \biggl( \sum_{n=1}^\infty \big( \frac{(-1)^{n-1}}{[n]} + \frac{Q^n}{[n]} \big) \frac{\beta_{-n}}{n} \biggr) \rangle.
\ee
\end{proposition}

Recall the following identities proved in \cite{Zho3}:
\be
W_{\mu^t\nu^t}(q) = (-1)^{|\mu|+|\nu|} W_{\mu,\nu}(q^{-1}),
\ee
and
\be
W_{\nu}(q) = (-1)^{|\nu|} q^{\kappa_\nu/2} W_\nu(q^{-1}).
\ee
So we have
\begin{multline} \label{eqn:Zb2}
\tilde{Z}^{(a)}(\lambda;t;\bp) \\
= \sum_{\mu, \nu, \eta} q^{(a+1) \kappa_\mu/2} (-1)^{|\mu|} W_{\mu\nu}(q^{-1}) Q^{|\nu|}q^{\kappa_\nu/2}  W_{\nu}(q^{-1})\cdot
\frac{\chi_{\mu}(\eta)}{z_{\eta}\sqrt{-1}^{l(\eta)}} p_{\eta}(\bx).
\end{multline}
This is essentially $Z^{-(a+1)}(-\lambda;t;\bp)$,
with an extra $(-1)^{|\mu|}$.
Hence by Proposition 3.1,
one can easily get:

\begin{proposition}
The normalized open string amplitude of the resolved conifold with one outer brane and framing $a \in \bZ$
is given by:
\be \label{eqn:ZbOperator1}
\hat{\tilde{Z}}^{(a)}(\lambda;t; \bp)
=  \langle \exp (\sum_{n=1}^\infty \frac{p_n(\bx)}{ni} \beta_n) q^{aK}
\exp \biggl( \sum_{n=1}^\infty \big( \frac{1}{[n]} - \frac{e^{-nt}}{[n]} \big) \frac{\beta_{-n}}{n} \biggr) \rangle.
\ee
\end{proposition}

\subsection{Proof of the full Mari\~no-Vafa Conjecture}
Write
\be
\hat{Z}^{(a)}(\lambda;t;\bp)
= \exp \sum_{\mu \in \cP^+} \sum_{g=0}^\infty  F^{(a)}_{g; k; \mu} \lambda^{2g-2} e^{-kt} p_{\mu}(\bx).
\ee
and
\be
\hat{\tilde{Z}}^{(a)}(\lambda;t;\bp)
= \exp \sum_{\mu \in \cP^+} \sum_{g=0}^\infty  \tilde{F}^{(a)}_{g; k; \mu} \lambda^{2g-2} e^{-kt} p_{\mu}(\bx),
\ee

\begin{theorem} \label{thm:MV}
For the resolved conifold with one outer brane,
the generating functions for the open Gromov-Witten invariants are given by:
\begin{multline}
\exp \sum_{\mu \in \cP^+} \sum_{g=0}^\infty  \sum_{k=0}^\infty i^{l(\mu)} F^{(a)}_{g; k; \mu} \lambda^{2g-2}
e^{(|\mu|/2-k)t} p_{\mu}(\bx) \\
= \sum_\mu s_{\mu}(\bx)  q^{(a+1)\kappa_\mu/2} \cdot \dim_q R_{\mu^t},
\end{multline}
\begin{multline}
\exp \sum_{\mu \in \cP^+} \sum_{g=0}^\infty  \sum_{k=0}^\infty i^{l(\mu)} \tilde{F}^{(a)}_{g; k; \mu} \lambda^{2g-2}
e^{(|\mu|/2-k)t} p_{\mu}(\bx) \\
= \sum_\mu s_{\mu}(\bx)  q^{a\kappa_\mu/2} \cdot \dim_q R_{\mu}.
\end{multline}
\end{theorem}

\begin{proof}
It is easy to see from (\ref{eqn:ZaOperator1}) and (\ref{eqn:ZbOperator1}) that
\begin{multline}
\exp \sum_{\mu \in \cP^+} \sum_{g=0}^\infty  \sum_{k=0}^\infty i^{l(\mu)} F^{(a)}_{g; k; \mu} \lambda^{2g-2} e^{-kt} p_{\mu}(\bx) \\
= \langle \exp (\sum_{n=1}^\infty \frac{p_n(\bx)}{n} \beta_n) q^{(a+1)K}
\exp \biggl( \sum_{n=1}^\infty \big( \frac{(-1)^{n-1}}{[n]} + \frac{Q^n}{[n]} \big) \frac{\beta_{-n}}{n} \biggr) \rangle,
\end{multline}
\begin{multline}
\exp \sum_{\mu \in \cP^+} \sum_{g=0}^\infty  \sum_{k=0}^\infty i^{l(\mu)} \tilde{F}^{(a)}_{g; k; \mu} \lambda^{2g-2} e^{-kt} p_{\mu}(\bx) \\
= \langle \exp (\sum_{n=1}^\infty \frac{p_n(\bx)}{n} \beta_n) q^{aK}
\exp \biggl( \sum_{n=1}^\infty \big( \frac{1}{[n]} - \frac{e^{-nt}}{[n]} \big) \frac{\beta_{-n}}{n} \biggr) \rangle,
\end{multline}
Now we have
\be
\exp (\sum_{n=1}^\infty \frac{p_n(\bx)}{n} \beta_n)
= \sum_\mu s_\mu(\bx) \langle \mu |.
\ee
Therefore by using (\ref{eqn:QDimT}) and (\ref{eqn:QDim2}),
the proof is completed.
\end{proof}

\end{document}